\documentclass[11pt,leqeq]{article}
\usepackage{amssymb,amsthm}
\usepackage{amsmath}
\usepackage{amscd}
\usepackage{xy}
\xyoption{all}
\usepackage[dvips]{graphicx}
\newtheorem{thm}{Theorem}
\newtheorem{lema}[thm]{Lemma}
\newtheorem{cor}[thm]{Corolary}

\newtheorem{defi}[thm]{Definition}
\newtheorem{exa}[thm]{Example}
\newcommand{\ob}{{\rm{Ob}}}

\newcommand{\inv}{\operatorname{inv}}
\newcommand{\val}{\operatorname{val}}
\begin{document}
\title{Feynman-Jackson integrals}
\author{Rafael D\'\i az\ \
and Eddy Pariguan} \maketitle
\begin{abstract}
We introduce perturbative Feynman integrals in the context of
$q$-calculus generalizing the Gaussian $q$-integrals introduced by
D\'\i az and Teruel. We provide analytic as well as combinatorial
interpretations for the Feynman-Jackson integrals.
\end{abstract}

\section{Introduction}
Feynman integrals are a main tool in high energy physics since
they provide a universal integral representation for the
correlation functions of any Lagrangian quantum field theory whose
associated quadratic form is non-degenerated. In some cases the
degenerated situation may be approached as well by including odd
variables as is usually done in the BRST-BV procedure. Despite its
power Feynman integrals still await a proper definition from a
rigorous mathematical point of view. The main difficulties in
understanding Feynman integrals are the following
\begin{enumerate}
\item {The output of a perturbative Feynman integral is a formal
power series in infinitely many variables, i.e., an element of
$\mathbb{C}[[g_1,...,g_n,..]].$ This fact goes against our
strongly held believe that the output of an integral should be a
number.}
\item {There is no guarantee that the formal series mentioned above
will be convergent, not even in an asymptotic sense. General
statements in this matter are missing.}

\item {Feynman integrals of greatest interest are performed over
spaces of infinite dimension. In this situation the coefficients
of the series in variables $\mathbb{C}[[g_1,...,g_n,..]]$ referred
above are given by finite dimensional integrals which might be
divergent. In this case additional care must be taken in order to
renormalize the values of these integrals. The renormalization
procedure, when applies, is done in two steps one of analytic
nature called regularization, and a further step of algebraic
nature which may be regarded as a fairly general form of the
inclusion-exclusion principle of combinatorics.}
\item {In process 1 to 3 above a number of choices must be made. No
general statements showing the unicity of the result are known.}
\end{enumerate}

{\noindent}Finite dimensional Feynman integrals are also of
interest for example in Matrix theory. They still present
difficulties 1 and 2 above but issues 3 and 4 become null. The
goal of this paper is to construct a $q$-analogue of Feynman
integrals which we call Feynman-Jackson integrals. We consider
only the simplest case of 1-dimensional integrals. Our approach is
to use the $q,k$-generalized gamma function and the
$q,k$-generalized Pochhamer symbol introduced in \cite{ED} and
\cite{CTT}.

{\noindent}The computation of a $1$-dimension Feynman integrals,
for example an integral of the form $\displaystyle{\int
e^{h(x)}dx}$ where $\displaystyle{h(x)=\frac{-x^2}{2}+
\sum_{j=1}^{\infty}h_j \frac{x^j}{j!}}$ is done in four steps

\begin{enumerate}
\item{ The integral is obtain perturbatively, meaning that the
integrand $h(x)$ should be replaced by a formal power series
$\displaystyle{\frac{-x^2}{2}+ \sum_{j=1}^{\infty}
g_jh_j\frac{x^j}{j!}\in \mathbb{C}[[g_1,...,g_n,..]]}$ where
$\{g_j\}_{j=1}^{\infty}$ is a countable set of independent
variables.} \item {One uses the key identity\ \
$e^{\frac{-x^2}{2}+ \sum_{j=1}^{\infty}
g_jh_j\frac{x^j}{j!}}=e^{-\frac{x^2}{2}}\ \ e^{\sum_{j=1}^{\infty}
g_jh_j\frac{x^j}{j!}}.$} \item {$e^{\sum_{j=1}^{\infty}g_jh_j
\frac{x^j}{j!}}$\ \ is expanded as a formal power series in
$\mathbb{C}[[g_1,...,g_n,..]]$ using the series expansion of
$e^x$. This step reduces the computation of the Feynman integrals
to computing a countable number of Gaussian integrals.} \item
{Compute the Gaussian integrals obtained in step $3$ which yield
as output an element in $\mathbb{C}[[g_1,...,g_n,..]]$.}
\end{enumerate}

{\noindent}Steps $1$,$3$ and $4$ can be carried out in
$q$-calculus without much difficulty . The subtle issue to be
tackled is the unpleasant fact that the identity
$e^{x+y}=e^{x}e^{y}$ does not hold in $q$-calculus.

{\noindent}This paper is organized as follows: in Section 2 after
a quick review of $q$-calculus we introduce Gauss-Jackson
integrals based on the definition of the function $\Gamma_{q,2}$
introduced in \cite{CTT}. In Section $3$ we introduce the
combinatorial tools that shall be needed to formulate our main
theorem. In Section $4$ we introduce the algebraic properties of
the $q$-exponential that will allow us to overcome the fact that
$E_q^{x+y}\neq E_q^{x}E_q^y$. In Section $5$ we shall enunciate
and prove our main result Theorem \ref{qanalo}
\begin{equation*}
\frac{1}{\Gamma_{q,2}(1)}\displaystyle{\int_{-\nu}^{\nu}
E_{q,2}^{\frac{-q^2x^2}{[2]_q}+\sum_{j=1}^{\infty} g_j h_j
\frac{x^j}{[j]_q!}}} d_qx=\sum_{\Lambda \in
\ob(\textbf{Graph}_q)/\sim} h_q(\Lambda)
\frac{\omega_q(\Lambda)}{aut_q(\Lambda)}
\end{equation*}
which gives a $q$-analogue of $1$-dimensional Feynman integrals.

\section{Gauss-Jackson integrals}\label{Basic}

In this section we review several definitions in $q$-calculus. The
reader may find more information in \cite{CTT}, \cite{Ch},
 \cite{K} and \cite{So} . We focus upon the $q,k$-generalizations
of the Pochhammer symbol, the gamma function and its integral
representations \cite{CTT}.

{\noindent} Let us fix $0<q<1$ and let
$f:\mathbb{R}\longrightarrow\mathbb{R}$ be any map. The
$q$-derivative $\partial_q(f)$ of $f$ is given by
$\displaystyle{\partial_q(f)=\frac{d_qf}{d_qx}=\frac{I_q(f)-
f}{(q-1)x}},$ \hspace{0.1cm} where
$I_q:\mathbb{R}\longrightarrow\mathbb{R}$ is given by
$I_q(f)(x)=f(qx)$ for all $x\in\mathbb{R},$ and $d_q(f)=I_q(f)-f$.

{\noindent} The definite Jackson integral (see \cite{Jac} and
\cite{JA}) of a map $f:[0,b] \longrightarrow \mathbb{R}$ is given
by
$$\int_{0}^{b}f(x)d_qx=(1-q)b\sum_{n=0}^{\infty}q^nf(q^nb).$$
{\noindent}The improper Jackson integral of a map $f:[0,\infty)
\longrightarrow
\mathbb{R}$ is given by
$$ \int_0^{\infty/a} f(x) d_qx =
(1-q)\sum_{n\in\mathbb{Z}}\frac{q^n}{a}f\left(\frac{q^n}{a}\right).$$
{\noindent}For all $t\in\mathbb{Z}^+$  the $q$-factorial is given
by $[t]_q!=[t]_q[t-1]_q \cdot \cdot \cdot [1]_q$, where
$\displaystyle{[t]_q=\frac{(1-q^t)}{(1-q)}}$ is the $q$-analogue
of a real number $t$. The $q$-factorial is an instance of the
$q,k$-generalized Pochhammer symbol which is given  by
$$\displaystyle{[t]_{n,k}= [t]_q[t+k]_q[t+2k]_q\dots[t+(n-1)k]_q =
\prod_{j=0}^{n-1}[t+jk]_q},\ \ \mbox{ for all
}\ t\in\mathbb{R}. $$

{\noindent}In this paper we shall mainly use the $q,2$-
generalized Pochhamer symbol evaluated at $t=1$, namely
\begin{equation}\label{pocha2}
[1]_{n,2}= [1]_q[3]_q[5]_q\dots[2n-1]_q =
\prod_{j=0}^{n-1}[1+2j]_q.
\end{equation}
{\noindent}We remark that $[1]_{n+1,2}= [2n+1]_q[1]_{n,2}.$ We
shall use the following notation. Let $x,y,t\in\mathbb{R}$ and $n
\in
\mathbb{Z^+}$ we set
$$\displaystyle{(x+y)_{q,2}^n := \prod_{j=0}^{n-1} (x+ q^{2j}y)}
\mbox{ and }
(1+x)_{q,2}^{t}:=\frac{(1+x)_{q,2}^{\infty}}{(1+q^{2t}x)_{q,2}^{\infty}},
\mbox{ where } {(1+x)_{q,2}^{\infty}:=\prod_{j=0}^{\infty} (1+
q^{2j}x)}$$

{\noindent}Recall that one can define two $q$-analogues of the
exponential function given as follows
\begin{eqnarray*}
E_{q,2}^{x} &=& \displaystyle{
\sum_{n=0}^{\infty}\frac{q^{n(n-1)}x^n}{[n]_{q^2}!}=(1+(1-q^2)x)_{q,2}^{\infty}}\\
e_{q,2}^x &=&
\displaystyle{\sum_{n=0}^{\infty}\frac{x^n}{[n]_{q^2}!}
=\frac{1}{(1-(1-q^2)x)_{q,2 }^{\infty}}}.
\end{eqnarray*}

{\noindent}The $q,2$-gamma function $\Gamma_{q,2}(t)$ is given by
the explicit formula
$\displaystyle{\Gamma_{q,2}(t)=\frac{{(1-q^2)_{q,2}^{{\frac{t}{2}}
-1}}}{{(1-q)^{\frac{t}{2}-1}}}}$ \quad for a real number $t>0$,
and has a representation in terms of $E_{q,2}^x$ given by the
following Jackson integral
\begin{equation}\label{gamatilde}
\displaystyle{\Gamma_{q,2}(t)=\int_{0}^{\left(\frac{[2]_q}{(1-q^2)}\right)^{\frac{1}{2}}}x^{t-1}E_{q,2}^{-\frac{q^2
x^2}{[2]_q}}d_qx}, \quad t>0.
\end{equation}
{\noindent} Similarly one can define $\gamma_{q,2}^{(a)}(t)$ for
$a>0$ using $e_{q,2}^x$ by the following Jackson integral

\begin{equation}\displaystyle{\gamma_{q,2}^{(a)}(t)=\int_{0}^{\infty/a(1-q^2)^{\frac{1}{2}}}
x^{t-1} e_{q,2}^{-\frac{x^2}{[2]_q}}}d_qx, \quad
t>0.\end{equation}

{\noindent}Both integral representations are related by
$\Gamma_{q,2}(t)=c(a,t)\gamma_{q,2}^{(a)}(t),$ where the function
$c(a,t)$ is given by
$$\displaystyle{c(a,t)=\frac{a^t[2]_q^{\frac{t}{2}}}{1+[2]_qa^{2}}\left(
1+ \frac{1}{[2]_qa^2}\right)_{q,2}^{\frac{t}{2}} \left( 1+[2]_qa^2
\right)_{q,2}^{1-\frac{t}{2}}}, \ \ \mbox{for}\  a>0\ \mbox{and}\  t\in \mathbb{R}.$$

{\noindent}We proceed to introduce two different $q$-analogues of
the Gaussian integral and given a Jackson integral representation
for each one. The Gaussian integrals  are related to each other by
the function $c(a,t)$ in a similar way as the integral
representations of the $q,2$-generalized gamma functions are
related to each other.

\begin{defi}\label{jgauss}
Let $\nu=\left(\frac{[2]_q}{(1-q^2)}\right)^{\frac{1}{2}}$ and
$\varepsilon^{(a)}=\infty/a(1-q^2)^{\frac{1}{2}}$, the
Gaussian-Jackson integrals are given by

\begin{equation*}\hspace{-1.3cm} G(t):=\frac{1}{2}\int_{-\nu}^{\nu}
x^{t-1}E_{q,2}^{-\frac{q^2 x^2}{[2]_q}}d_qx=\frac{1}{2}
\int_{0}^{\nu} x^{t-1}E_{q,2}^{-\frac{q^2
x^2}{[2]_q}}d_qx+ \frac{1}{2}\int_{-\nu}^{0}
x^{t-1}E_{q,2}^{-\frac{q^2 x^2}{[2]_q}}d_qx,  \quad t>0.
\end{equation*}

\begin{equation*}\hspace{-0.6cm}
G^{(a)}(t):=\frac{1}{2}\int_{-\varepsilon^{(a)}}^{\varepsilon^{(a)}}
x^{t-1}
e_{q,2}^{-\frac{x^2}{[2]_q}}d_qx=\frac{1}{2}\int_{0}^{\varepsilon^{(a)}}
x^{t-1} e_{q,2}^{-\frac{x^2}{[2]_q}}d_qx +\frac{1}{2}
\int_{-\varepsilon^{(a)}}^{0} x^{t-1}
e_{q,2}^{-\frac{x^2}{[2]_q}}d_qx, \quad t>0.
\end{equation*}

\end{defi}
{\noindent} Notice that if $t-1$ is an odd integer both integrals
in Definition \ref{jgauss} are zero because then $x^{t-1}$ is an
odd function while $ E_{q,2}^{-\frac{q^2 x^2}{[2]_q}}$ and
$e_{q,2}^{-\frac{x^2}{[2]_q}}$ are even functions.

\section{Combinatorial interpretation of $[1]_{n,2}$}\label{comb}

In this section we introduce the combinatorial tools that will be
needed in order to describe $q$-analogue of $1$-dimensional
Feynman integrals. The interested reader may consult \cite{AND},
\cite{RA}, \cite{HY} for further information.

\begin{defi}
A partition of $a\in \mathbb{Z}^+$ is a finite sequence of
positive integers $a_1, a_2,\dots,a_r$ such that
$\displaystyle{\sum_{i=1}^{r}a_i=a}.$ For $a,d\in \mathbb{Z}^+$,
$p_d(a)$ denotes the number of partitions of $a$ into less than
$d$ parts.

\end{defi}

\begin{defi}
Let $n\in \mathbb{Z}^+$ and $a_1, a_2,\dots,a_n$  be a partition
of $a$. The $q$-multinomial coefficient is given by
$$\left[%
\begin{array}{c}
a_1 + a_2 +\dots + a_n  \\
  a_1, a_2, \dots ,a_n \\
\end{array}%
\right]_q=\frac{[a_1+a_2+\dots a_n]_q!}{[a_1]_q![a_2]_q!
\dots [a_n]_q!} . $$
\end{defi}

{\noindent} Denote by $[[n]]$ the set $\{1,\dots,n\}$ ordered in
the natural way. $|X| $ denotes the cardinality of set $X$ and
$S_X$ denotes the group of permutations on $X$.
\begin{defi}\label{inver}
Let $a_1,\dots,a_n$ be a partition of $a$. We denote by
$S(a_1,\dots,a_n)$ the set of all maps $f :\ [[a]]
\longrightarrow [[n]]$ such that $|f^{-1}(i)|=a_i$ for all $i \in
[[n]]$. We set $\inv(f):= |\{(i,j) \in [[a]]\times[[a]] : i<j
\mbox{ and } \ f (i)
>f(j)\}|.$
\end{defi}
{\noindent}The following result is proved using induction.
\begin{thm}\label{thecorrs}
$$\left[%
\begin{array}{c}
a_1 + a_2 +\dots + a_n  \\
  a_1, a_2,\dots,a_n \\
\end{array}%
\right]_q=\sum_{f\in S(a_1,\dots,a_n)}q^{\inv (f)}.$$

\end{thm}

{\noindent} Notice that this result  implies that $$[n]_q!=\left[%
\begin{array}{c}
1 + 1 +\dots + 1  \\
  1, 1, \dots ,1 \\
\end{array}%
\right]_q=\sum_{f\in S(1,\dots,1)}q^{\inv (f)}=\sum_{f \in
S_{[[n]]}}q^{\inv (f)}.$$
\begin{defi}
A paring $\alpha$ on a totally ordered set $R$ of cardinality $2n$
is a sequence $\alpha=\{ \ (a_i,b_i) \ \}_{i=1}^{n} \in (R^2)^n$
such that
\begin{enumerate}
\item {$a_1< a_2 < \dots < a_n$.}
\item{ $a_i<b_i,\quad
i=1,\dots,n.$}
\item{ $\displaystyle{R=\bigsqcup_{i=1}^{n} \{a_i,b_i\}.}$}
\end{enumerate}

\end{defi}

{\noindent}We denote by $P(R)$ the set of pairings on $R$.

\begin{defi}
For $\alpha \in P([[2n]])$ we set
\begin{enumerate}
\item {$((a_i,b_i))=\{j\in [[2n]]: a_i<j<b_i\}$ for all $(a_i,b_i)\in
\alpha$.}
\item {$P_i(\alpha)=\{b_j: 1 \leq j < i \}$.}
\item {$\displaystyle{w(\alpha)=\prod_{i=1}^{n}q^{|((a_i,b_i)) \setminus
P_i(\alpha)|}=q^{\sum_{i=1}^{n}|((a_i,b_i)) \setminus
P_i(\alpha)|}}.$ We call $w(\alpha)$ the weight of $\alpha$.}
\end{enumerate}
\end{defi}

\begin{figure}[h!]
\begin{center}
\includegraphics{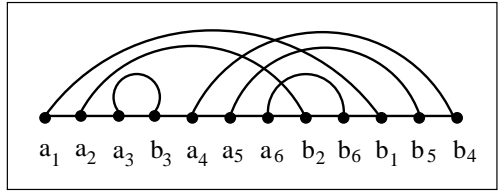}
\end{center}
\caption{Example of a pairing $\alpha$ on $[[12]]$}
\label{fig:weights}
\end{figure}

\begin{exa}
Let $\alpha$ be the pairing on $[[12]]$  shown in Figure
\ref{fig:weights}. The weight $w(\alpha)$ can be computed as follows

\[
\begin{array}{lll}
q^{|((a_1,b_1)) \setminus P_1(\alpha)|}=q^8, & q^{|((a_2,b_2))
\setminus P_2(\alpha)|}=q^5, & q^{|((a_3,b_3))
\setminus P_3(\alpha)|}=q^0, \\
\mbox{ }&\mbox{ }&\mbox{ }\\
q^{|((a_4,b_4)) \setminus P_4(\alpha)|}=q^{6-2}, & q^{|((a_5,b_5))
\setminus P_5(\alpha)|}=q^{4-2}, & q^{|((a_6,b_6)) \setminus
P_6(\alpha)|}=q^{1-1}.
\end{array}
\]

Hence $w(\alpha)=q^{19}.$
\end{exa}

\newpage
\begin{thm}\label{pochateor}
Given $n\in\mathbb{N}$ the following identity holds
\begin{equation}\label{intpocha} [1]_{n,2}=\sum_{\alpha \in
P([[2n]])}w(\alpha).
\end{equation}

\end{thm}
\begin{proof}
We use induction on $n$. For $n=1$, we have $[1]_{1,2}=1$. Suppose
identity (\ref{intpocha}) holds for $n$, we prove it for $n+1$ as
follows

\begin{eqnarray*}
\sum_{\alpha \in P([[2n+2]])} w(\alpha)&=&\sum_{\alpha \in
P([[2n+2]])} w(\alpha-\{(a_1,b_1)\})q^{|((a_1,b_1))|}\\&=&\sum_{2
\leq b_1 \leq 2n+2}q^{b_1-2} \sum_{\beta \in P([[2n+2]]\setminus
\{(a_1,b_1)\} )} w(\beta) \\&=& \sum_{2 \leq b_1 \leq 2n+2} q^{b_1
- 2} \sum_{\beta \in P([[2n]])} w(\beta)\\&=& \sum_{2\leq b_1 \leq
2n+2} q^{b_1 - 2} \quad
[1]_{n,2}\\&=&[2n+1]_q[1]_{n,2}=[1]_{n+1,2}.
\end{eqnarray*}

\end{proof}

Notices that as $q\longrightarrow 1$ we recover the well known
identity
$$(2n-1)(2n-3)\dots 1=|\{\mbox{pairings on}\ [[2n]]\}|.$$
\begin{exa}By definition  $[1]_{2,2}=[1]_q[3]_q=[3]_q.$ Consider the
pairings of a four elements ordered set. Figure \ref{fig:Pochha
graphs} shows that there are $3$ such pairings and that the sum of
their weights is $1+q+q^2$ as it should.
\end{exa}

\begin{figure}[h!]
\begin{center}
\includegraphics{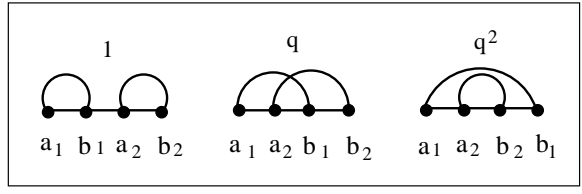}
\end{center}
\caption{Combinatorial meaning of $[1]_{2,2}$} \label{fig:Pochha
graphs}
\end{figure}

\section{Algebraic properties of the $q$-exponentials}

The $q$-exponential maps $e_q^x$ and $E_q^x$ are good
$q$-analogues of the exponential map $e^x$ since they satisfy
$\partial_qe_q^x=e_q^x$ , $e_q^0$=1 and
$\displaystyle{\lim_{q\longrightarrow 1}e_q^x=e^x},$ and
$\partial_qE_q^x=E_q^{qx}$, $E_q^0=1$ and
$\displaystyle{\lim_{q\longrightarrow 1}E_q^x=e^x}$. From a
differential point of view $e_q^x$ is the right $q$-analogue of
$e^x$. However both $e_q^x$ and $E_q^x$ lack the fundamental
algebraic property of the exponential, namely that
$e^x:(\mathbb{R},+) \longrightarrow (\mathbb{R},\cdot) $ is a
group homomorphism. Indeed one checks that $e_q^{x+y}\neq e_q^x
e_q^y$ and also that $E_q^{x+y}\neq E_q^x E_q^y.$ Nevertheless we
still have the remarkable identity $(e_q^x)^{-1}=E_q^{-x}.$

{\noindent} A possible algebraic solution to this problem is to
assume that $yx=qxy$. Using this relation one verifies that
$e_q^{x+y}=e_q^x e_q^y$ and $E_q^{x+y}= E_q^x E_q^y.$ However we
still have to deal with the fact that $e_q^{x+y}\neq e_q^x e_q^y$
and $E_q^{x+y}\neq E_q^x E_q^y$ for commuting variables $x,y\in
\mathbb{R}$. Theorems \ref{separa} and \ref{separa2} below provide
tools that allow us to overcome this obstacle in the process of
computing Feynman integrals, as discussed in the Introduction.

\begin{thm}\label{separa}
$\displaystyle{E_{q,2}^{x+y}=E_{q,2}^{x}\left(\sum_{c,d\geq 0}
{\lambda}_{c,d}x^{c}y^{d} \right)}$, where
$\displaystyle{{\lambda}_{c,d}=\sum_{k=0}^{c}\frac{(-1)^{c-k}{{d +k}
\choose k}q^{(d +k)(d +k - 1)}}{[d + k]_{q^2}![c - k]_{q^2}!}.}$
\end{thm}
\begin{proof}
\vspace{-0.5cm}
\begin{eqnarray*}
E_{q,2}^{x+y}e_{q,2}^{-x}&=&\left({\sum_{n=0}^{\infty}}\frac{q^{n(n-1)}(x+y)^n}{[n]_{q^2}!}\right)
\left(\sum_{m=0}^{\infty}\frac{(-1)^m x^m}{[m]_{q^2}!}\right)\\
&=&\sum_{n,m,k\leq n}\frac{(-1)^m {n \choose k}q^{n(n-1)}}{[n]_{q
^2}![m]_{q^2}!}x^{m+k}y^{n-k}
\end{eqnarray*}

{\noindent} Making the change $c=m+k \mbox{ and } d=n-k$, we get
\begin{equation}
E_{q,2}^{x+y}e_{q,2}^{-x}= \sum_{c,d\geq
0}\left(\sum_{k=0}^{c}\frac{(-1)^{c-k}{{d +k} \choose k}q^{(d +k)(d
+k - 1)}}{[d + k]_{q^2}![c - k]_{q^2}!}\right)x^{c}y^{d}.
\end{equation}
\end{proof}

\begin{thm}\label{separa2}
$\displaystyle{e_{q,2}^{x+y}=e_{q,2}^{x}\left(\sum_{c,d\geq 0}
{\kappa}_{c,d}x^{c}y^{d} \right)}$, where
$\displaystyle{{\kappa}_{c,d}=\sum_{k=0}^{c}\frac{(-1)^{c-k}{{d
+k}
\choose k}q^{(c-k)(c-k - 1)}}{[d + k]_{q^2}![c - k]_{q^2}!}}.$
\end{thm}
\begin{proof}
\begin{eqnarray*}
e_{q,2}^{x+y}E_{q,2}^{-x}&=&\left(\sum_{n=0}^{\infty}\frac{(x+y)^n}{[n]_{q^2}!}\right)
\left(\sum_{m=0}^{\infty}\frac{(-1)^m q^{m(m-1)} x^m}{[m]_{q,2}!}\right)\\
&=&\sum_{n,m,k\leq n}\frac{(-1)^m {n \choose
k}q^{m(m-1)}}{[n]_{q^2}![m]_{q^2}!}x^{m+k}y^{n-k}
\end{eqnarray*}

{\noindent} Fixing $c=m+k \mbox{ and } d=n-k$, we have
\begin{equation}
e_{q,2}^{x+y}E_{q,2}^{-x}=\sum_{c,d\geq
0}\left(\sum_{k=0}^{c}\frac{(-1)^{c-k}{{d +k} \choose k}q^{(c-k)(c-k
- 1)}}{[d + k]_{q^2}![c - k]_{q^2}!}x^{c}y^{d}\right).
\end{equation}
\end{proof}

\begin{lema}\label{limi} For $c,d\in\mathbb{N}$, ${\displaystyle {\lim_{q\rightarrow 1}\lambda_{c,d}=\frac{1}{d!}
\delta_{c,0}}}$.
\end{lema}
\begin{proof}
$${\displaystyle {\lim_{q\rightarrow 1}\lambda_{c,d}= \sum_{k=0}^{c} (-1)^{c-k}
\frac{{d+k\choose k}}{(d+k)!(c-k)!}=\frac{1}{d!c!}\sum_{k=0}^{c} (-1)^{c-k} {c\choose k}=\frac{1}{d!}
\delta_{c,0}}}.$$
\end{proof}
\section{Feynman-Jackson integrals}

We denote by  \textbf{Graph} the category whose objects
$\ob(\textbf{Graph})$ are graphs. Recall that a graph $\Lambda$ is
triple $(V,E,b)$ where $V$ and $E$ are finite sets, called the set
of vertices and the set of edges respectively, and $b$ is a map
that assigns to each edge $e \in E$ a subset of $V$ a cardinality
one or two. To each graph we associate a map $\val: V
\longrightarrow \mathbb{N}$ defined by $\val(s)=|\{e: s\in
b(e)\}|$. All graphs considered in this paper are such that $\val
(s)\geq 1$ for all $s \in V $.
{\noindent} Morphisms in
\textbf{Graph} from $\Lambda_1$ to $\Lambda_2$ are pairs
$(\varphi_V,\varphi_E)$ such that

\begin{enumerate}
\item {$\varphi_V: V(\Lambda_1) \longrightarrow V(\Lambda_2)$.}
\item {$\varphi_E: E(\Lambda_1) \longrightarrow E(\Lambda_2)$.}
\item {$b({\Lambda_2})(\varphi_E(e))=\varphi_V(b({\Lambda_1})(e)), \mbox{ for
all } e\in E(\Lambda_1)$.}
\end{enumerate}
The essence of $1$-dimensional Feynman integrals, see \cite{ET},
may be summarized in the following identity
\begin{equation}\label{fey1}
\frac{1}{\sqrt{2\pi}}\int e^{\frac{-x^2}{2}+
\sum_{j=1}^{\infty}g_jh_j\frac{x^j}{j!}}=\sum_{\Lambda \in
\ob(\textbf{Graph})/\sim} h(\Lambda)
\frac{\omega(\Lambda)}{aut(\Lambda)}.
\end{equation}
In identity (\ref{fey1}) the following notation is used
\begin{enumerate}
\item{$\ob(\textbf{Graph})/\sim$ denotes the set of isomorphisms classes of graphs.}
\item {$h(\Lambda)=\prod_{s \in V}h_{\val(s)}$.}
\item {$\omega(\Lambda)= \prod_{s \in V} g_{\val(s)}$.}
\item {$aut(\Lambda)=|Aut(\Lambda)|$
where $Aut(\Lambda)$ denotes the set of isomorphisms from graph
$\Lambda$ into itself, for all $\Lambda \in \ob(\textbf{Graph})$}.

\end{enumerate}
{\noindent} Theorem \ref{qanalo} below provides a $q$-analogue of
identity (\ref{fey1}). We first prove the following
\begin{thm}\label{separ}
\begin{equation*}
\frac{1}{\Gamma_{q,2}(1)}\displaystyle{\int_{-\nu}^{\nu}
E_{q,2}^{\frac{-q^2x^2}{[2]_q}+\sum_{j=1}^{\infty} g_j h_j
\frac{x^j}{[j]_q!}}} d_qx=\sum_{m=0}^{\infty}\chi_m q^m
\end{equation*}
where

$$\chi_m= \sum_{\alpha
,c,d,j,k,l,f}\frac{(-1)^{k}g_lh_l}{[2]_{q}^{c}[2j]_q![d+k]_{q^2}![c-k]_{q^2}!}{{d+k}
\choose k}. $$ The sum above runs over all $c,d,j,k \in
\mathbb{Z}^+$, such that $k\leq c$, $l\in p_d(2j)$, $f \in
S(l_1,\dots,l_d)$, $\alpha \in P([[2c+2j]])$ and
$\displaystyle{\sum_{i=1}^{c+j}|((a_i ,b_i))\setminus
P_i(\alpha)|+ \inv(f) + (d+k)(d+k-1)+2c=m}.$

\end{thm}
\begin{proof}
Making the changes $\displaystyle{x \longrightarrow
\frac{-q^2x^2}{[2]_q}}$ and $\displaystyle{y \longrightarrow
\sum_{j=1}^{\infty}\frac{g_jh_jx^j}{[j]_q!}}$ in Theorem
\ref{separa}, we get
\begin{eqnarray*}
& &E_{q,2}^{\frac{-q^2x^2}{[2]_q}+\sum_{j=1}^{\infty}g_jh_j
\frac{x^j}{[j]_q!}}=E_{q,2}^{\frac{-q^2x^2}{[2]_q}}\sum_{c,d=0}^{\infty}\left(\frac{\lambda_{c,d}(-1)^cq^{2c}x^{2c}}{[2]_q^c}\left(\sum_{j=1}^{\infty}g_jh_j
\frac{x^j}{[j]_q!}\right)^d
\right)\\&=&E_{q,2}^{\frac{-q^2x^2}{[2]_q}}\sum_{c,d=0}^{\infty}\left(\frac{\lambda_{c,d}(-1)^cq^{2c}x^{2c}}{[2]_q^c}\left(\sum_{j=1}^{\infty}\left(\sum_{l
\in p_d(j)}\frac{[j]_{q}!g_{l_1}\dots g_{l_d}h_{l_1}\dots
h_{l_d}}{[l_1]_q![l_2]_q!\dots  [l_d]_q!}\right)
\frac{x^j}{[j]_q!}\right)\right)\\&=&
E_{q,2}^{\frac{-q^2x^2}{[2]_q}}\sum_{c,d,j,l}\frac{\lambda_{c,d}(-1)^cq^{2c}g_{l_1}\dots
g_{l_d}h_{l_1}\dots h_{l_d}}{[2]_q^c}\left[%
\begin{array}{c}
j  \\
  l_1, \dots ,l_{d} \\
\end{array}%
\right]_q x^{2c+j}.
\end{eqnarray*}
Using the expression given in Theorem \ref{separa} for
$\lambda_{c,d}$ and using the convention that $g_l = g_{l_1}\dots
g_{l_d} $ and $h_l = h_{l_1}\dots  h_{l_d} $ for $l\in p_d(j)$ we
get

\begin{equation}\label{preint}
E_{q,2}^{\frac{-q^2x^2}{[2]_q}+\sum_{j=1}^{\infty}g_jh_j
\frac{x^j}{[j]_q!}}=E_{q,2}^{\frac{-q^2x^2}{[2]_q}}\sum_{c,d,j,k,l}
\frac{(-1)^{2c-k} g_l h_{l}q^{(d+k)(d+k-1)+2c}}{[2]_q^c
[j]_q![d+k]_{q^2}![c-k]_{q^2}!}{{d+k \choose k}}
\left[%
\begin{array}{c}
j  \\
  l_1,\dots,l_{d} \\
\end{array}%
\right]_q x^{2c+j}.
\end{equation}
{\noindent} Multiply by $\frac{1}{\Gamma_{q,2}(1)}$ and integrate
both sides of the equation (\ref{preint}) from $-\nu$ to $\nu$
(which cancels out all terms with $j$ odd), one gets for $l \in
p_d(2j)$

\begin{eqnarray}\label{primera}
&&
\frac{1}{\Gamma_{q,2}(1)}\int_{-\nu}^{\nu}E_{q,2}^{\frac{-q^2x^2}{[2]_q}+
\sum_{j=1}^{\infty}g_jh_j
\frac{x^j}{[j]_q!}} d_qx\\ \hspace{-0.5cm}\label{secun}&=&
\hspace{-0.6cm}\small \sum_{c,d,j,k,l} \frac{(-1)^{2c-k}
g_lh_{l}q^{(d+k)(d+k-1)+2c}{d+k \choose k}}{[2]_q^c
[2j]_q![d+k]_{q^2}![c-k]_{q^2}!}{}
\left[%
\begin{array}{c}
2j  \\
  l_1,\dots ,l_{d} \\
\end{array}%
\right]_q
\frac{1}{\Gamma_{q,2}(1)}\int_{-\nu}^{\nu}E_{q,2}^{\frac{-q^2x^2}{[2]_q}}x^{2c+2j}d_qx\\
\label{ter}&=&\hspace{-0.6cm} \sum_{c,d,j,k,l} \frac{(-1)^{k}
g_lh_{l}q^{(d+k)(d+k-1)+2c}{d+k \choose k}}{[2]_q^c
[2j]_q![d+k]_{q^2}![c-k]_{q^2}!}{}
\left[%
\begin{array}{c}
2j  \\
  l_1,  \dots ,l_{d} \\
\end{array}%
\right]_q [1]_{c+j,2}.
\end{eqnarray}
{\noindent} Notice that equation (\ref{ter}) is obtained from
equation (\ref{secun}) using Definition \ref{jgauss}. Using
Theorem
\ref{pochateor} and (\ref{inver}) in the right-hand side of
(\ref{ter}) one obtains

\begin{equation}\label{odimf}
\sum_{\alpha,c,d,j,k,l,f}\frac{(-1)^{2c-k}g_lh_l{{d+k} \choose
k}}{[2]_{q}^{c}[2j]_q![d+k]_{q^2}![c-k]_{q^2}!}\displaystyle{q^{\sum_{i=1}^{c+j}|((a_i
,b_i))\setminus P_i(\alpha)|+ \inv(f) + (d+k)(d+k-1)+2c}}.
\end{equation}
Which yields the desired result.
\end{proof}
Using Lemma \ref{limi} one notices that the limit as $q$ goes to
$1$ of (\ref{odimf}) is
$${\displaystyle{\sum \frac{g_{l}h_{l}}{(2j)!d!} {2j\choose {l_1,\dots,l_d}} |{\mbox{pairings on}\ [[2j]]}|}}$$
which is well known to be equivalent to formula (\ref{fey1})

\begin{defi}
We denote by  $\textbf{Graph}_q$ the category whose objects
$\ob(\textbf{Graph}_q)$ are planar q-graphs $(V,E,b,f)$ such that

\begin{enumerate}
\item {$V= \{ \bullet \} \sqcup V^{1} \sqcup V^{2}$ where
$V^{1}=\{\otimes_1,\dots,\otimes_{|V_1|}\}$ and
$V^{2}=\{\circ_1,\dots,\circ_{|V_2|}\}$. } \item{ $E=E^1 \sqcup
E^2 \sqcup E^3$.} \item{$b$ is a map that assigns to each edge
$e\in E$ a subset of $V$ a cardinality two.} \item{ Set
$F_{\circ}=\{(\circ_i,e): i\in [[|V^2|]]\ \mbox{and}\ \bullet
\notin b(e)\}$. We require that $|F_{\circ}|$ be even.
$f:F_{\circ}\longrightarrow [[|V^2|]]$ is any map.} \item
{$|b^{-1}(\{\otimes_i,\bullet \})|\in \{0,1\}$ for all $i\in
[[|V^1|]]$ and $|b^{-1}(\circ_i,\bullet)|=1$ for all $i\in
[[|V^2|]].$ If $|b^{-1}(\{ \otimes_i,\bullet \} )|=1$ then
$|b^{-1}(\{ \otimes_j,\bullet \} )|=1$ for all $j\geq i$; and
$\bullet\in b(e)$ for any $e\in E^{3}$. } \item {$\val
(\otimes_i)\in \{2,3\}$. If $\val(\otimes_i)=3$ then
$\val(\otimes_j)=3$ for all $i\geq j$, and $|E^2|\leq |V^{1}|$.}
\end{enumerate}

\end{defi}

\begin{figure}[h]
\begin{center}
\includegraphics{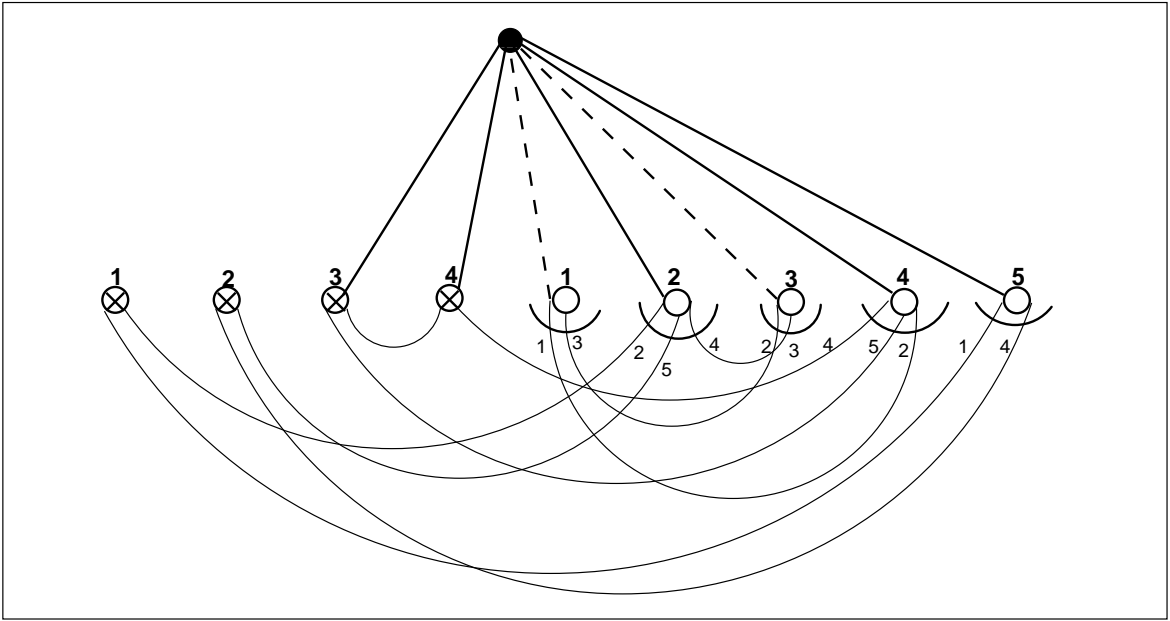}
\caption{Feynman $q$-diagram} \label{fig:quantumgraph}
\end{center}
\end{figure}

{\noindent} Morphisms in $\textbf{Graph}_q$ are defined in the
obvious way. Figure \ref{fig:quantumgraph} shows an example of a
planar $q$-graph with $n=4$ and $m=5.$ Edges in $E^1$ $(E^2, E^3$)
are depicted by dark (dotted, regular) lines, respectively. The
map $f$ can be read off the numbering of half-edges in $E^3$
attached to vertices $\{\circ_1,\dots,\circ_m\}$.

{\noindent}Notice that associated to any graph
$\Lambda\in\ob(\textbf{Graph}_q)$ there exists a pairing $\alpha$
on the naturally ordered set $\{(v,e): v\in V^{1}\sqcup V^{2} \
\mbox{and}\ \bullet\notin b(e)\}$. Similarly, associated to any
graph there is a map $\widehat{f}:[[|F_{\circ}|]]\longrightarrow
[[|V^{2}|]]$ which is constructed from $f$ and the natural
ordering on $F_{\circ}$. For $\Lambda\in \ob({\textbf{Graph}_q})$
we set
\begin{enumerate}
\item{$h_q={\displaystyle \prod_{i=1}^{|V^{2}|} h_{\val(\circ_i)-1} }$.}
\item{$\omega_q ={\displaystyle(-1)^{|E^{2}|}
q^{2|V^{1}|+{ |V^{2}|+|E^{2}|\choose 2
}}\omega(\alpha)\inv(\widehat{f})\prod_{i=1}^{|V^{2}|}
g_{\val(\circ_i)-1} }$.}
\item{$aut_q(\Lambda)=[2]_q^{n}\ [|F_{\circ}|]_q!\
[|V^{2}|+|E^{2}|]_{q^{2}}!\
[|V^{1}|-|E^{2}|]_{q^{2}}!$}
\end{enumerate}

{\noindent} Using the notion of planar $q$-graphs introduced above
Theorem \ref{separ} may be rewritten as follows

\begin{thm}\label{qanalo}
\begin{equation*}
\frac{1}{\Gamma_{q,2}(1)}\displaystyle{\int_{-\nu}^{\nu}
E_{q,2}^{\frac{-q^2x^2}{[2]_q}+\sum_{j=1}^{\infty} g_j h_j
\frac{x^j}{[j]_q!}}} d_qx=\sum_{\Lambda \in
\ob(\textbf{Graph}_q)/\sim} h_q(\Lambda)
\frac{\omega_q(\Lambda)}{aut_q(\Lambda)}
\end{equation*}

\end{thm}

{\noindent} Setting $h_j=1$ in Theorem   \ref{qanalo} one gets
\begin{cor}

\begin{equation*}
\frac{1}{\Gamma_{q,2}(1)}\displaystyle{\int_{-\nu}^{\nu}
E_{q,2}^{\frac{-q^2x^2}{[2]_q}+\sum_{j=1}^{\infty}  g_j
\frac{x^j}{[j]_q!}}} d_qx=\sum \frac{{\displaystyle(-1)^{|E^{2}|}
q^{2|V^{1}|+{ |V^{2}|+|E^{2}|\choose 2
}}\omega(\alpha)\inv(\widehat{f})\prod_{i=1}^{|V^{2}|}
g_{\val(\circ_i)-1} }}{[2]_q^{n}\ [|F_{\circ}|]_q!\
[|V^{2}|+|E^{2}|]_{q^{2}}!\ [|V^{1}|-|E^{2}|]_{q^{2}}!}
\end{equation*}
where the sum runs over all ${\Lambda \in
\ob(\textbf{Graph}_q)/\sim}$.
\end{cor}

\section*{Acknowledgment}

Many thanks to Carolina Teruel.
\bibliographystyle{amsplain}
\bibliography{jackson-feynman}

$$\begin{array}{c}
  \mbox{Rafael D\'\i az. Universidad Central de Venezuela (UCV).} \ \  \mbox{\texttt{rdiaz@euler.ciens.ucv.ve}} \\
  \mbox{Eddy Pariguan. Universidad Central de Venezuela (UCV).} \ \
\mbox{\texttt{eddyp@euler.ciens.ucv.ve}} \\
\end{array}$$

\end{document}